\documentclass[12pt]{article}

\usepackage[top=0.75in]{geometry}
\usepackage{amsmath,amsthm,amssymb}
\usepackage{enumerate, graphicx}
\usepackage{xcolor}

\theoremstyle{definition}

\newtheorem{conj}{Conjecture}

\begin{document}

\title{A Quantitative Study on Average Number of Spins of Two-Player Dreidel}
\author{Thotsaporn ``Aek'' Thanatipanonda\\
Science Division,
Mahidol University International College\\
Nakornpathom, Thailand \\
{\tt thotsaporn@gmail.com}
}
\date{June 28, 2019}

\maketitle
\begin{abstract}
We give a high precision approximation of the average number of spins in a simplified version of a two-player version of the game Dreidel. We also make a conjecture on the average number of spins of the full version of the game.
\end{abstract}

\section{Introduction}
In 2006, Robinson and Vijay \cite{RV} showed that the average number of spins of the Dreidel with $k \geq 2$ players where each player starts with $n$ nuts is $\mathcal{O}(n^2)$. In this paper, we give a quantitative analysis of the average number of spins of a two-player simplified version of this game, where players start with $a$ and $b$ nuts, respectively.

\textbf{Motivation: Gambler's Ruin:}

We will first motivate our work with a well known statistical concept, gambler's ruin, which is interpreted in terms of a game that starts with a gambler who has zero chips.  For each play, s/he could gain one chip or lose one chip with equal probability $p=q=1/2$. The game continues until the gambler gains $M$ chips or loses $N$ chips. We are interested in the expected number of plays until the game ends.

\noindent Let $G(a)$ be the expected number of plays until the game ends when the gambler initially has $a$ chips. The system of recurrence relations that $G$ satisfies is as follows: 
\[ G(M)=  G(-N) =0,  \]
\[   G(a)=1+\dfrac{1}{2}G(a+1)+\dfrac{1}{2}G(a-1), -N <a < M.    \]

\noindent The solution to this system is well known, 
\begin{equation} \label{Gam}
 G(a)= (N+a)\cdot(M-a) ,  \;\  \;\   -N \leq a \leq M,
\end{equation}

\noindent i.e. \cite{Ross}. This beautiful result is a cornerstone of probability models 
and is the source of inspiration for our research problem.


\section{Analysis of Simplified Dreidel}
\vspace{0.5cm}
\begin{center}
\includegraphics[scale=0.52]{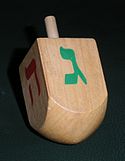}
\end{center}

A Dreidel is a four-sided spinning top, and the game associated with it is played during the Jewish holiday of Hanukkah. Each side of the Dreidel bears a letter of the Hebrew alphabet: gimel, hay, nun and shin. This is an $n$-player pot game where  player $i$ begins with $a_i$ nuts, contributing or taking away nuts from the pot as the game progresses. At the start of the game, everyone donates one nut to the pot and takes their turn to spin the Dreidel. The player who spins the Dreidel takes a certain action depending on how the Dreidel lands:
\begin{itemize}
\item Gimel: Player takes the whole pot, after which everyone donates one nut to the pot and then the next person spins.
\item Hay: Player takes the smaller half of the pot and the next person spins.
\item Nun: Player takes and gives nothing and the next person spins.
\item Shin: Player gives one nut to the pot and the next person spins.
\end{itemize}
The game ends when one person has all of the nuts in their possession.

\noindent We consider the simplest non-trivial version where two players spin with only two outcomes: Gimel (player takes the whole pot) and Shin (player gives one to the pot), and we answer the question of how long the game will last if the two players start with $a$ and $b$ nuts, respectively.

Let $D(a,p,b)$ denote the average number of spins until the game ends where the first player has $a$ nuts, the second player has $b$ nuts, and $p$ nuts are in the pot. We can derive the following recurrences for the simplified game:
\begin{align}
  D(a,p,b) &= 1+\dfrac{1}{2}\left[  D(b-1,2,a+p-1) + D(b,p+1,a-1)   \right],  \label{EqDr} \\
  D(0,p,b) &= D(a,p,0) =0 , \;\ \;\  \text{ for any } a \geq 0,\, p \geq 0,\, b\geq 0. \nonumber   
\end{align} 

\textbf{Example.} Some values of $D(a,p,b)$ when $a+b+p=6$ are
\[   D(1,1,4) = \dfrac{33}{16},  \;\ D(1,2,3) = \dfrac{5}{2},  \;\ D(1,3,2) = \dfrac{9}{4},  \;\ D(1,4,1) = 1,  \;\  \]
\[   D(2,1,3) = \dfrac{57}{16},  \;\ D(2,2,2) = 3,  \;\ D(2,3,1) = \dfrac{3}{2},  \;\  \]
\[   D(3,1,2) = \dfrac{15}{4},  \;\ D(3,2,1) = \dfrac{17}{8},  \;\ D(4,1,1) = \dfrac{9}{4}.  \]

\textbf{Experimental Math in Action:}

We are interested in the average number of spins until the game ends, \[T(a,b) := D(a,2,b).\] It would be too much to ask for an exact formula for $T(a,b)$ because we can see that the recurrence \eqref{EqDr} is quite complicated. The best we could hope for is a good approximation of it. The first thing that comes to mind is to compare $T(a,b)$ with the result of gambler's ruin (i.e., see \eqref{Gam}). It is reasonably safe to guess that $T(a,b) \leq ab$ for $a,b \geq 0$. 

\noindent After experimenting with different values of $D(a,p,b)$, we see that the sequence \[ \Delta_a := D(a+1,p,b)-D(a,p,b)\] quickly converges to a constant for each $p,b.$ A similar thing happens for \[ \Delta_b := D(a,p,b+1)-D(a,p,b).\] This implies that \[ \Delta_a^{2} \approx 0 , \;\ \;\  \Delta_b^{2} \approx 0.  \] In other words, 
$$D(a+2,p,b)-2D(a+1,p,b)+D(a,p,b) \approx 0$$ or
$$T(a+2,b)-2T(a+1,b)+T(a,b) \approx 0.$$ 

As a consequence, we learn that $T(a,b)$ is almost linear in both $a$ and $b$. $T(a,b)$ must be of the form \begin{equation} \label{conj} T(a,b) = c_3\cdot ab + c_2\cdot a + c_1\cdot b +c_0 +\epsilon_{a,b}, \end{equation} 
where $\epsilon_{a,b}$ is some small error in which we shall see that it is exponentially small in $a$ and $b$, i.e. the leading terms are $1/4^a$ and $1/4^b.$  
This is a fairly nice observation, but how could we prove it rigorously?

\textbf{Symbolic Computation in Action:}

Since we are only interested in $T(a,b)$ (in the situation where $p=2$), we will rewrite the system of equations \eqref{EqDr} recursively in terms of $T(a,b)$ only. 

\textbf{Example:}
\begin{align*}
T(1,3) &= 1+\dfrac{1}{2}\left[ T(2,2)+0   \right], \\
T(2,2) &= 1+\dfrac{1}{2}\left[ T(1,3)+D(2,3,1) \right] 
= 1+\dfrac{1}{2}\left[ T(1,3)+ 1+\dfrac{1}{2}[ 0+1] \right] \\
&= 1+ \dfrac{1}{2}+\dfrac{1}{4}+\dfrac{1}{2}T(1,3), \\
T(3,1) &= 1+\dfrac{1}{2}\left[ T(0,4)+D(1,3,2) \right]
 = 1+\dfrac{1}{2}\left[ T(0,4)+ 1+\dfrac{1}{2}[ T(1,3)+0] \right] \\
 &= 1+ \dfrac{1}{2}+\dfrac{1}{2}T(0,4)+\dfrac{1}{4}T(1,3)
 = 1+ \dfrac{1}{2}+\dfrac{1}{4}T(1,3).
\end{align*}

The process appears to get out of hand pretty quickly. But we can sort and check these relations with the actual values using the computer. As a result, we obtain
\begin{equation}   \label{GenT}    \tag{Key}
T(a,b) = \sum_{i=0}^{min(2a-2,2b-1)}\left(\dfrac{1}{2}\right)^i +
 \sum_{i=1}^{min(a,b)}  \dfrac{T(b-i,a+i)}{2^{2i-1}}+
\sum_{i=2}^{min(a,b+1)}  \dfrac{T(a-i,b+i)}{2^{2i-2}} .
\end{equation}          

The first sum comes from both players landing on shin (pay one) until one of them runs out of nuts. The second/third sum comes from the first/second player landing on gimel (takes the whole pot) respectively for the first time.

\textbf{The Guess and Check Method:}

So far so good! Equation \eqref{GenT} is the key! We simply plug in the conjectured equation \eqref{conj} into \eqref{GenT}. We consider two different cases according to the upper limits of the sums:
case 1: $a \geq b+1$,
case 2: $a \leq b.$

Indeed, with the help of Maple, the conjectured polynomial for $T(a,b)$ fits \eqref{GenT}, where the terms with $1/4^a$ and $1/4^b$ have been ignored. All of these cases produce the same solution where $c_3= 12/19$ and $c_1=c_2+2/19$ and there are no restrictions on $c_2$ and $c_0$. 

\textbf{The Approximation of $T(a,b)$:} 

The two free variables $c_2$ and $c_0$ could be approximated using the least squares method. In principle, we should work on the two cases separately. However, doing it all at once gives a very nice approximation, so we choose to do it this way. We use the least squares method over the values of $T(a,b)$  where $30 \leq a,b \leq 60$ to obtain 
\begin{align*}
T(a,b) \approx \widetilde T(a,b) =  \dfrac{12}{19}ab+c_2a +\left(c_2+\dfrac{2}{19}\right)b+c_0,   
\end{align*}
in which  $$c_2 = -0.304636562751640396971893222635\dots$$
and
$$c_0= 2.13102617218341081870452144156\dots .$$ 

This approximation is astonishing, for example, 
$|T(100,100)-\widetilde T(100,100)| < 10^{-12}.$

\textbf{The approximation of $D(a,p,b)$!:} 

We don't only have a good approximation for every starting position
but also a good approximation for every position!! 
For any numbers of nuts $a$ and $b$, we notice that the sequence 
\begin{equation} \label{clue1}
 h_{a,b}(p) := D(a,p+1,b)-D(a,p,b)
 \end{equation}
converges to a constant. 
Moreover these constants are linear in $a$ and $b$ i.e.
\begin{equation} \label{clue2}
 h_{a,b}(p) = s_2\cdot a + s_1 \cdot b + s_0.
 \end{equation}
 By combining these two ideas i.e. \eqref{clue1} and \eqref{clue2},
 it is reasonable to set up the solution form of $D(a,p,b)$ as 
\begin{equation} \label{MyD}   \tag{GOOD GUESS}
   D(a,p,b)  = T(a,b) + (p-2)(s_2\cdot a + s_1 \cdot b + s_0)  + \epsilon_{a,p,b}.
 \end{equation}
This is a nice idea! But how can we prove it? We make use of the original recurrence relation
\[  D(a,p,b) = 1+ \dfrac{1}{2}D(b-1,2,a+p-1)+\dfrac{1}{2}D(b,p+1,a-1).  \]
Plug in \eqref{MyD} on both sides and equate the coefficients of 
$ab, a ,b$ and the constant term. The conjectured equation fits. We obtain that 
$s_2 = \dfrac{4}{19}, s_1 = \dfrac{8}{19}$ and
$s_0 = c_2-\dfrac{18}{19}.$ The final solution is
\begin{align*}
D(a,p,b)
&\approx \dfrac{12}{19}ab+\left(c_2+\dfrac{4p}{19}-\dfrac{8}{19}\right)a 
+\left(c_2+\dfrac{8p}{19}-\dfrac{14}{19}\right)b+c_0 + \left(p-2\right)\left(c_2-\dfrac{18}{19}\right),
\end{align*}
where $c_2$ and $c_0$ were defined earlier. 
We stress that the approximation works better for large $a$ and $b$.

\textbf{The Main Conjecture:} 

\noindent After successfully finding a good solution to the simplified Dreidel game, it is natural to try to see whether the same method works for the full version (that is, where all four outcomes can occur). Let us denote the average number of spins for the full version by $Dr(a,p,b).$ We notice that each of the sequences
\[\Delta_a := Dr(a+1,p,b)-Dr(a,p,b)\] 
and 
\[  \Delta_b := Dr(a,p,b+1)-Dr(a,p,b).\]
converge to constants as before. Therefore, the same technique should work once we figure out relations similar to \eqref{GenT}. For now, we only conjecture how long the Dreidel game lasts and denote
$Q(a,b) := Dr(a,2,b).$ \\

\begin{conj}[Dreidel conjecture]
\begin{align*}
 Q(a,b) \approx \widetilde {Q}(a,b) &= 
2.21814151862618181904832628843 \cdot ab \\
& -1.09709667033405910669478639669 \cdot a \\
& -0.447079544643588135688652268182 \cdot b \\
& +2.83880783734231869675987135868 .
\end{align*}
\end{conj}

We obtain $\widetilde {Q}(a,b)$ by the least-squares method on $15 \leq a,b \leq 25.$ This function agrees with Zeilberger \cite{Z1}, if we let $a=b=nuts-1.$ The approximation also comes with incredible accuracy, for example 
$|Q(35,22)-\widetilde {Q}(35,22)| < 10^{-10}.$

For a practical perspective, 
we may assume that it takes 10 seconds per play,
if both players start with 10 nuts, 
an average game will last 28.10 minutes.  
And if both players start with 15 nuts, 
an average game will last 69.33 minutes.

\end{document}